\documentclass{article}
\usepackage[margin=1.6in]{geometry}
\usepackage{amsmath}
\usepackage{amssymb}
\usepackage{amsthm}
\usepackage{graphicx}
\usepackage{color}
\usepackage{algorithm}
\usepackage{url}
\usepackage{algpseudocode}
\usepackage[font=footnotesize]{caption}
\algtext*{Comment}
\algnewcommand{\Comment}[1]{\hfill// #1}

\newcommand{\R}{{\mathbb R}}

\theoremstyle{plain}
\newtheorem{theorem}{Theorem}[section]

\theoremstyle{definition}
\newtheorem{definition}[theorem]{Definition}

\theoremstyle{remark}
\newtheorem{remark}[theorem]{Remark}

\definecolor{BLUE}{RGB}{0,0,255}

\usepackage{fancyhdr}
\begin{document}
\title{3D Printing of Invariant Manifolds in Dynamical Systems} 
\author{Patrick R. Bishop, Summer Chenoweth, Emmanuel Fleurantin,\\ Alonso Ogueda-Oliva, Evelyn Sander, and Julia Seay} 
\date{}
\maketitle

\vspace{-2.5em}
\begin{center}
 {\small Department of Mathematical Sciences, George Mason University, Fairfax, VA 22030}
\end{center}

\begin{abstract}
Invariant manifolds are one of the key features that organize the dynamics of a differential equation. 
We introduce a novel approach to visualizing and studying invariant manifolds by using 3D printing technology, combining advanced computational techniques 
with modern 3D printing processes to transform mathematical abstractions into tangible models. Our work addresses the challenges of translating complex manifolds into printable meshes, showcasing results for the following systems of differential equations: the Lorenz system, the Arneodo-Coullet-Tresser system, and the Langford system. By bridging abstract mathematics and physical reality, this approach promises new tools for research and education in nonlinear dynamics. We conclude with practical guidelines for reproducing and extending our results, emphasizing the potential of 3D-printed manifolds to enhance understanding and exploration in dynamical systems theory.\end{abstract}


\begin{figure}[!htb]
    \centering
    \includegraphics[width=.5\textwidth]{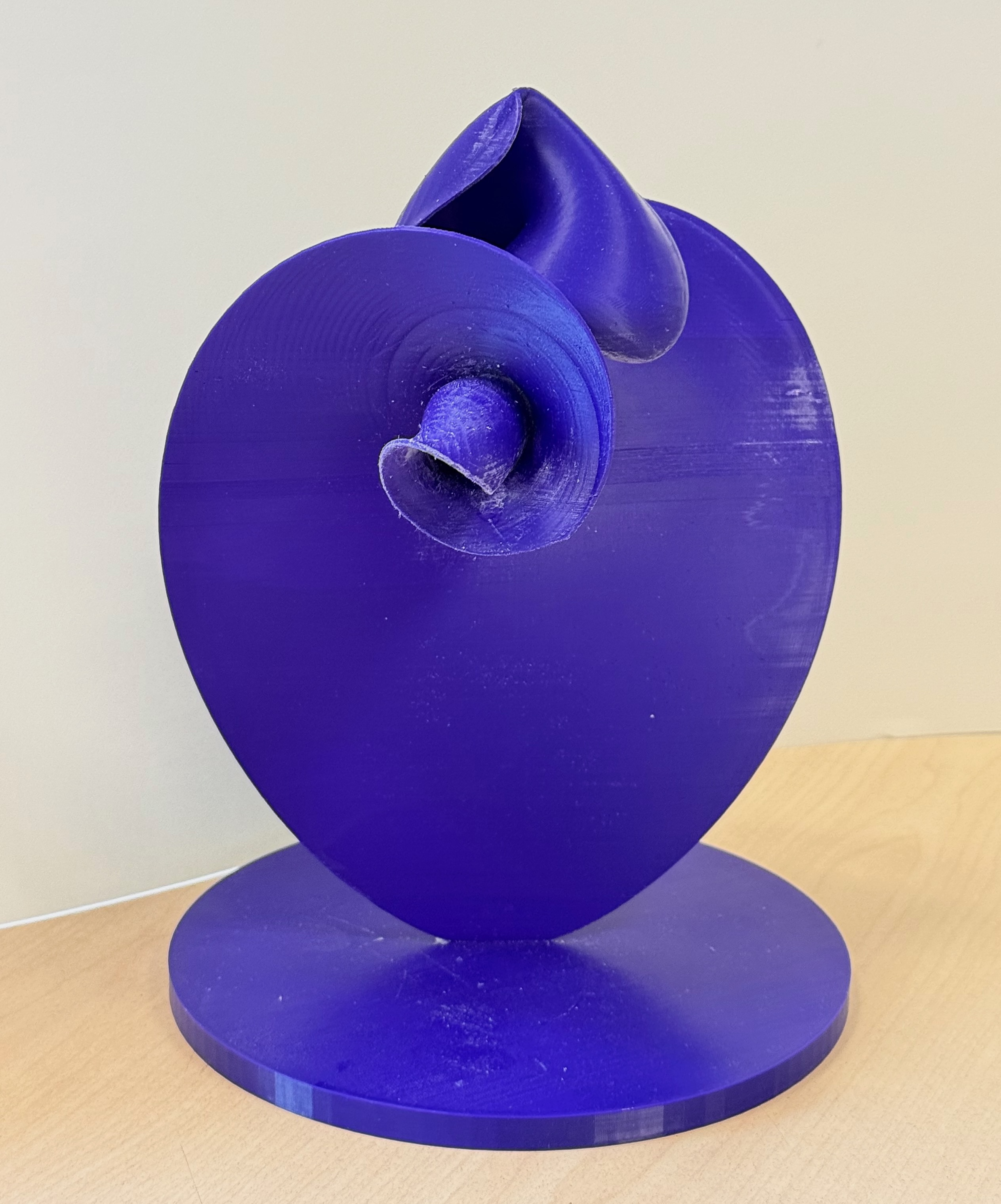}
    \caption{A 3D printed stable manifold of the origin for the Lorenz system ~\eqref{eq:lorenz}.}
    \label{fig:phase_portrait_lorenz}
\end{figure}

\begin{figure}[!htb]
\tabcolsep=0pt
\begin{center}
\begin{tabular}{cc}
\includegraphics[height=5cm, width=6cm, keepaspectratio]{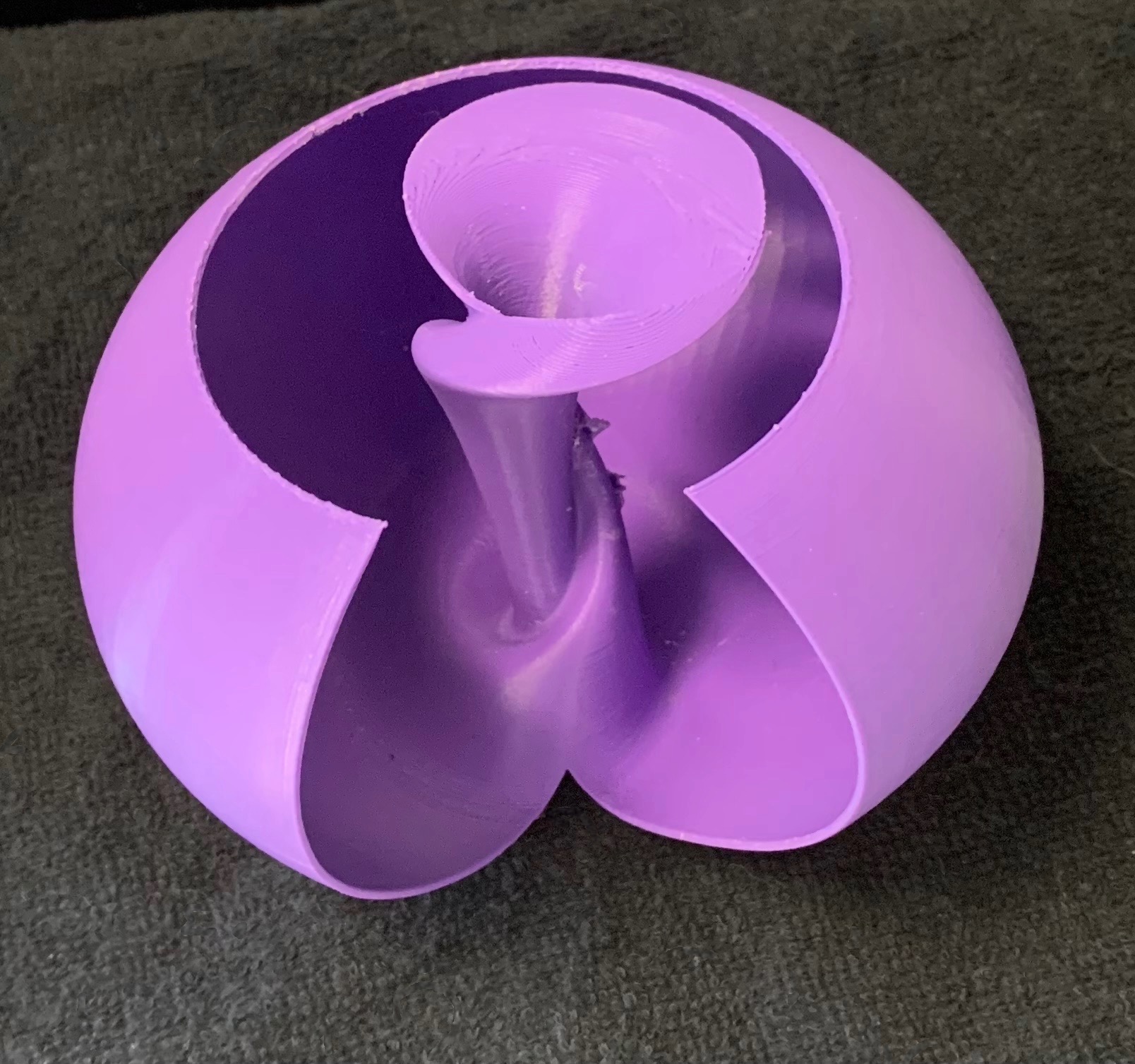} & 
\hspace{0.5em}\includegraphics[height=5cm, width=6cm, keepaspectratio]{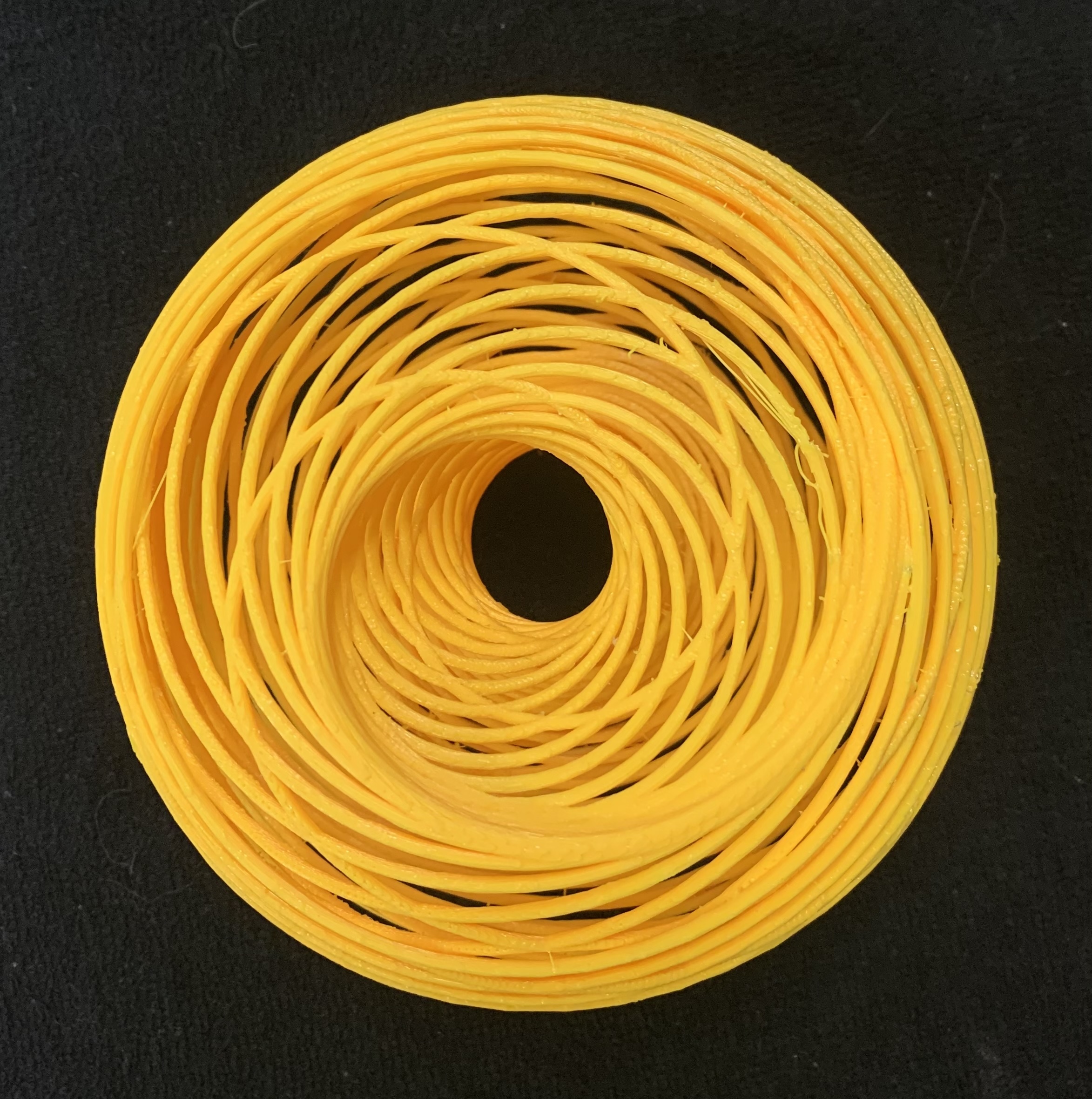} \\ 
\end{tabular}
\vspace{0.5cm} 
\begin{tabular}{c}
\includegraphics[height=5cm, width=6cm, keepaspectratio]{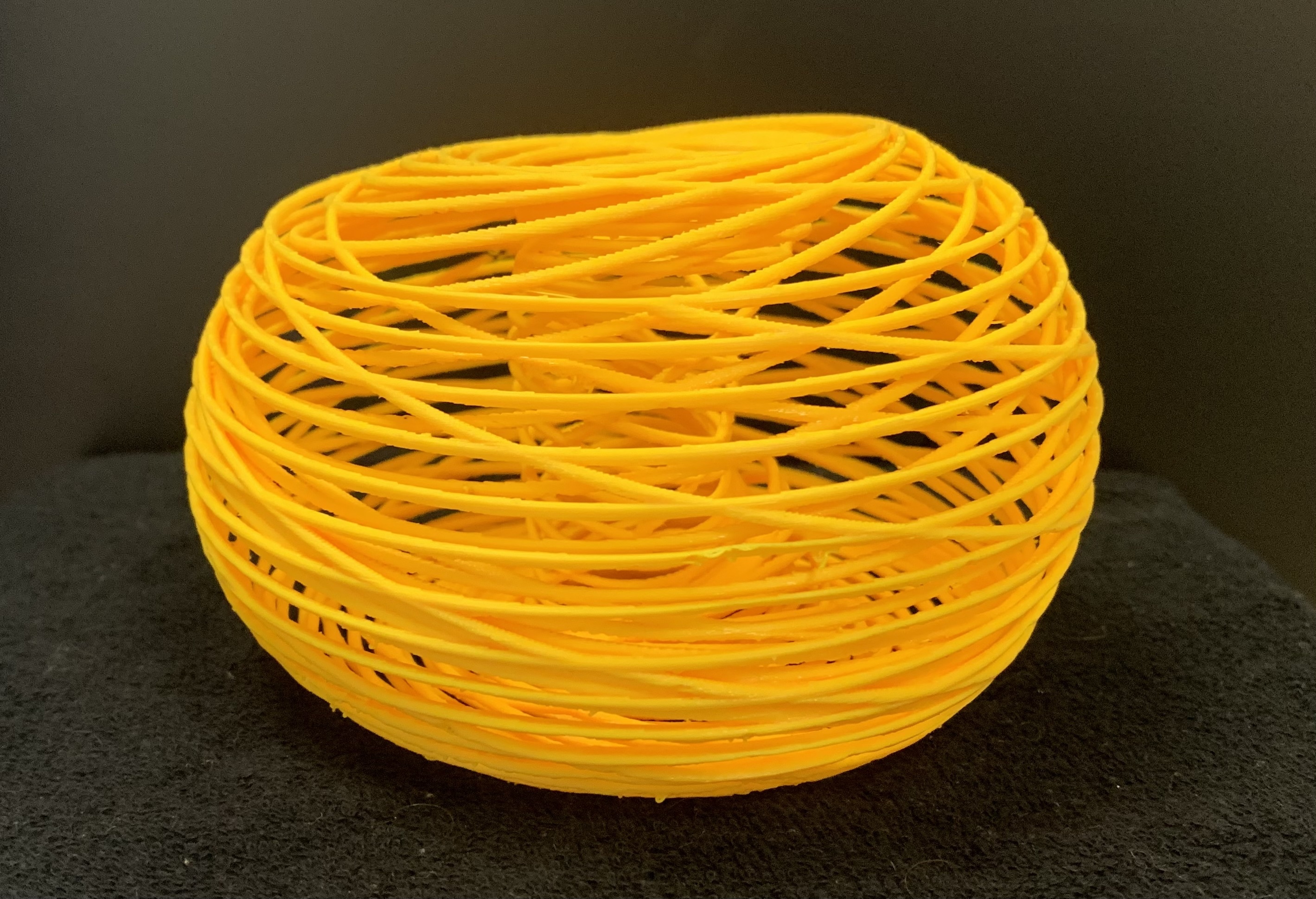} \\
\end{tabular}
\end{center}
\caption{This picture shows 3D printed unstable manifolds of equilibria in the Langford system~\eqref{eq:langford}, displaying structural changes as the bifurcation parameter $\alpha$ varies. Left: The unstable manifold of $p_1 \approx 1.94$ for $\alpha=0.95$ exhibits a complex structure with multiple lobes. Creating a cut-out allows the viewer to see the internal structure of the manifold. Top-right and bottom: Two views of the unstable manifold of $p_1 \approx 1.84$ for $\alpha=0.806$, revealing its spiral structure from above (top-right) and its three-dimensional profile (bottom). Printing using a mesh structure is an alternative to the cut-out for revealing the manifold's intricate details.}
\label{fig:phase_portrait_langford}
\end{figure}

\section{Introduction}

Invariant manifolds play a crucial role in understanding the behavior of dynamical systems, providing geometric structures that organize the state space and reveal fundamental aspects of a system's long-term dynamics. These invariant sets, particularly stable and unstable manifolds, organize the dynamical behavior, such as serving as the separation between different basins of attraction. Through their intersections and foliations, the manifolds determine the qualitative structure of the flow. While computational methods for studying these structures have advanced significantly, visualizing and intuitively grasping their complex geometries remains challenging. 

\begin{figure}[!htb]
\tabcolsep=0pt
\begin{center}
\begin{tabular}{cc}
\includegraphics[height=6cm, width=7cm, keepaspectratio]{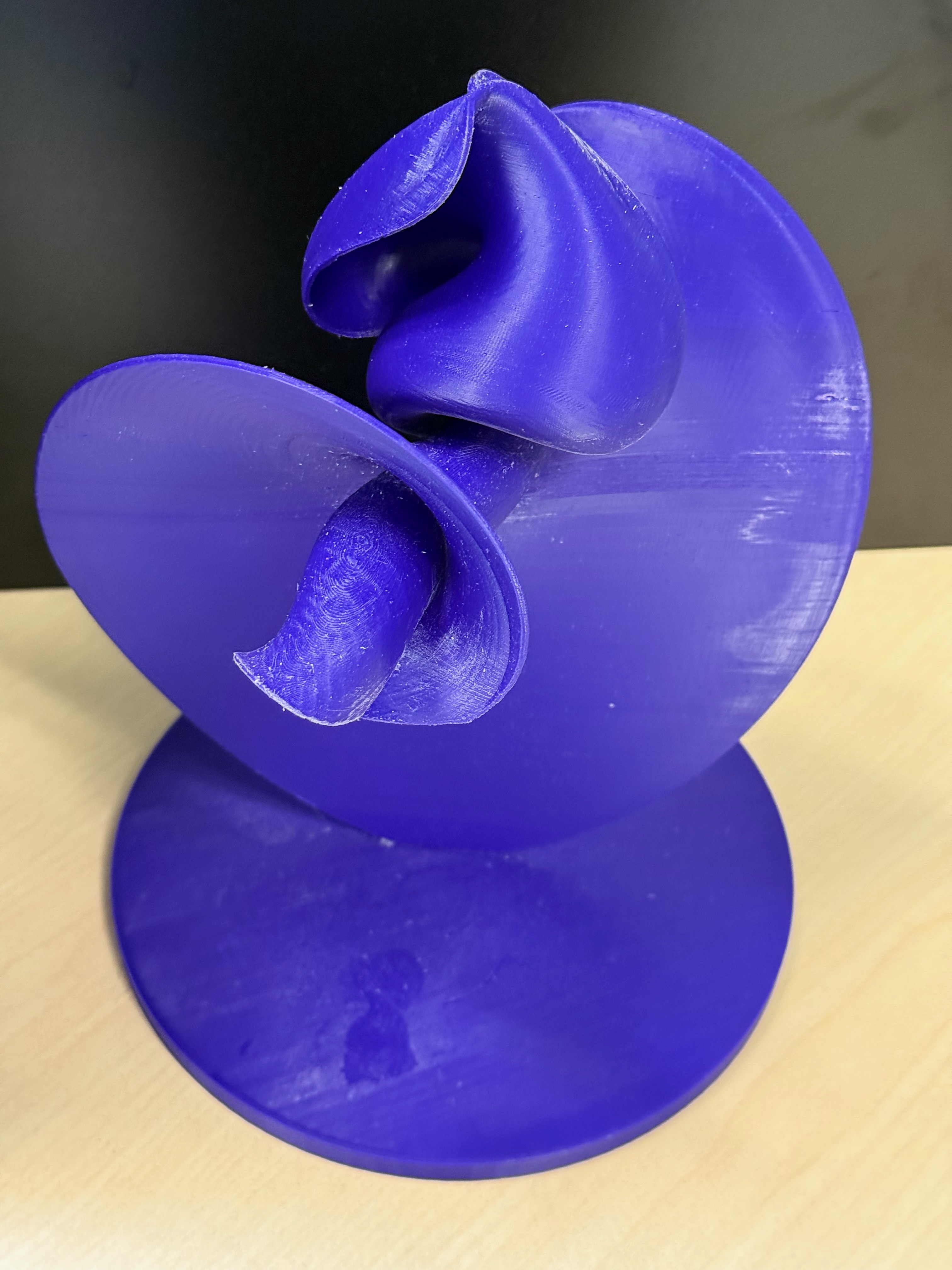} & 
\hspace{0.5em}\includegraphics[height=6cm, width=7cm, keepaspectratio]{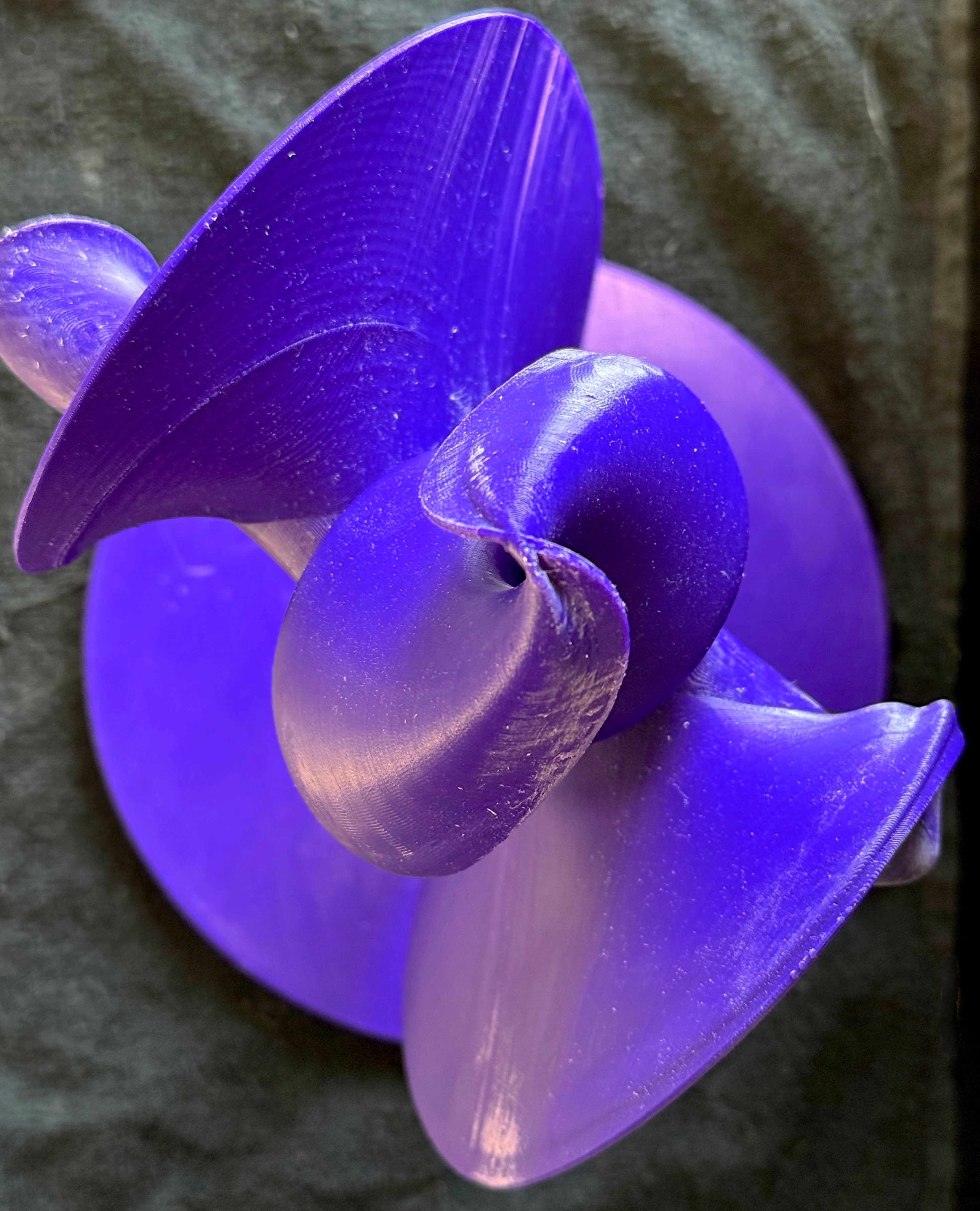} \\ 
\end{tabular}
\begin{tabular}{c}
\includegraphics[height=9cm, width=10cm, keepaspectratio]{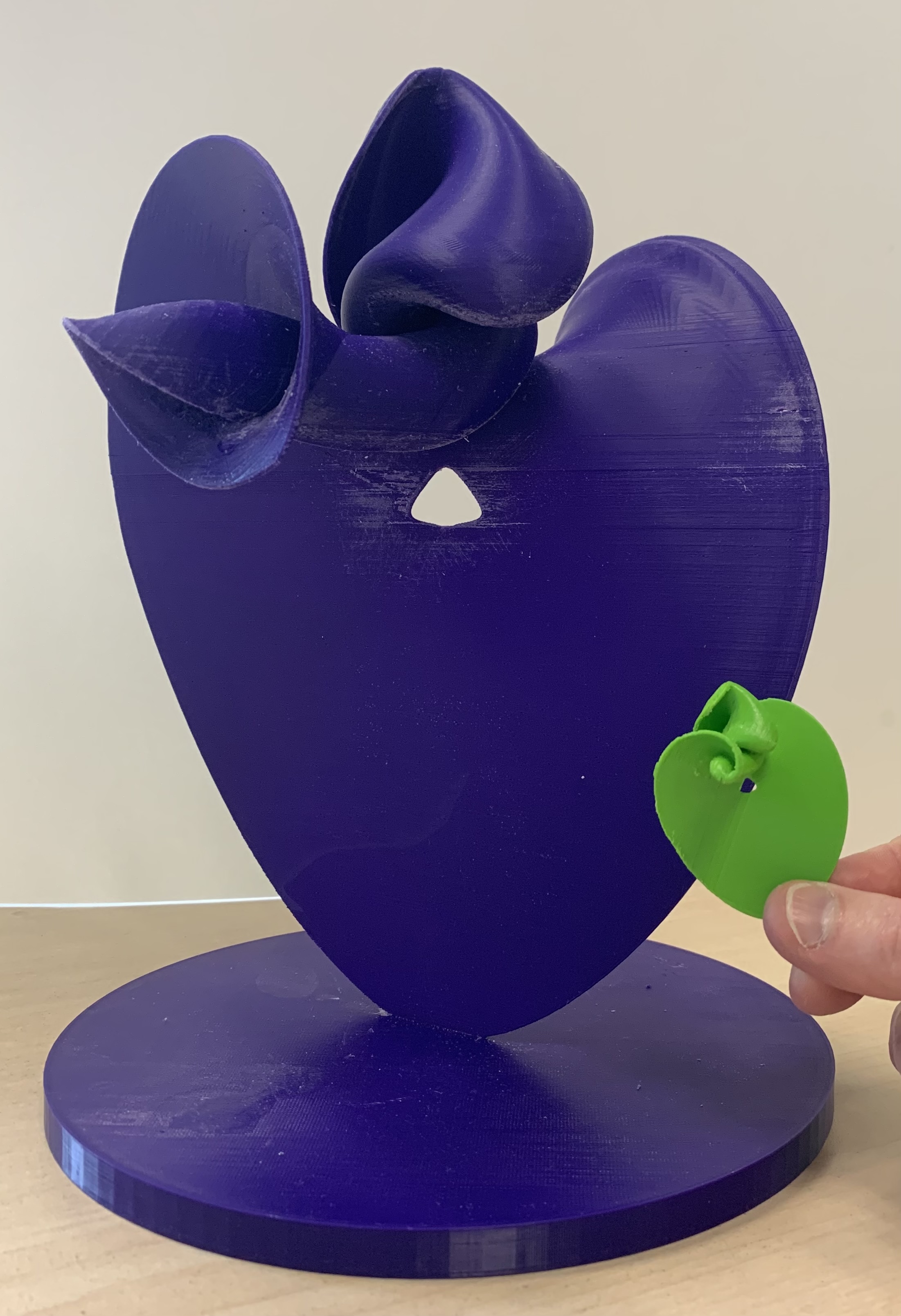} \\
\end{tabular}
\end{center}
\caption{Different perspectives of the 3D-printed stable manifold of the Lorenz system ~\eqref{eq:lorenz} at the origin. Top left and right: Two views of the larger model showing the characteristic spiral structure and complex geometry of the stable manifold. Bottom: Comparison with a smaller version, demonstrating that the same data can be used to print in a variety of sizes. 
The small version is held in hand for size reference. \label{fig:3dprint_manifold}}
\end{figure}

\begin{figure}[!htb]
\centering
\includegraphics[width=0.7\textwidth]{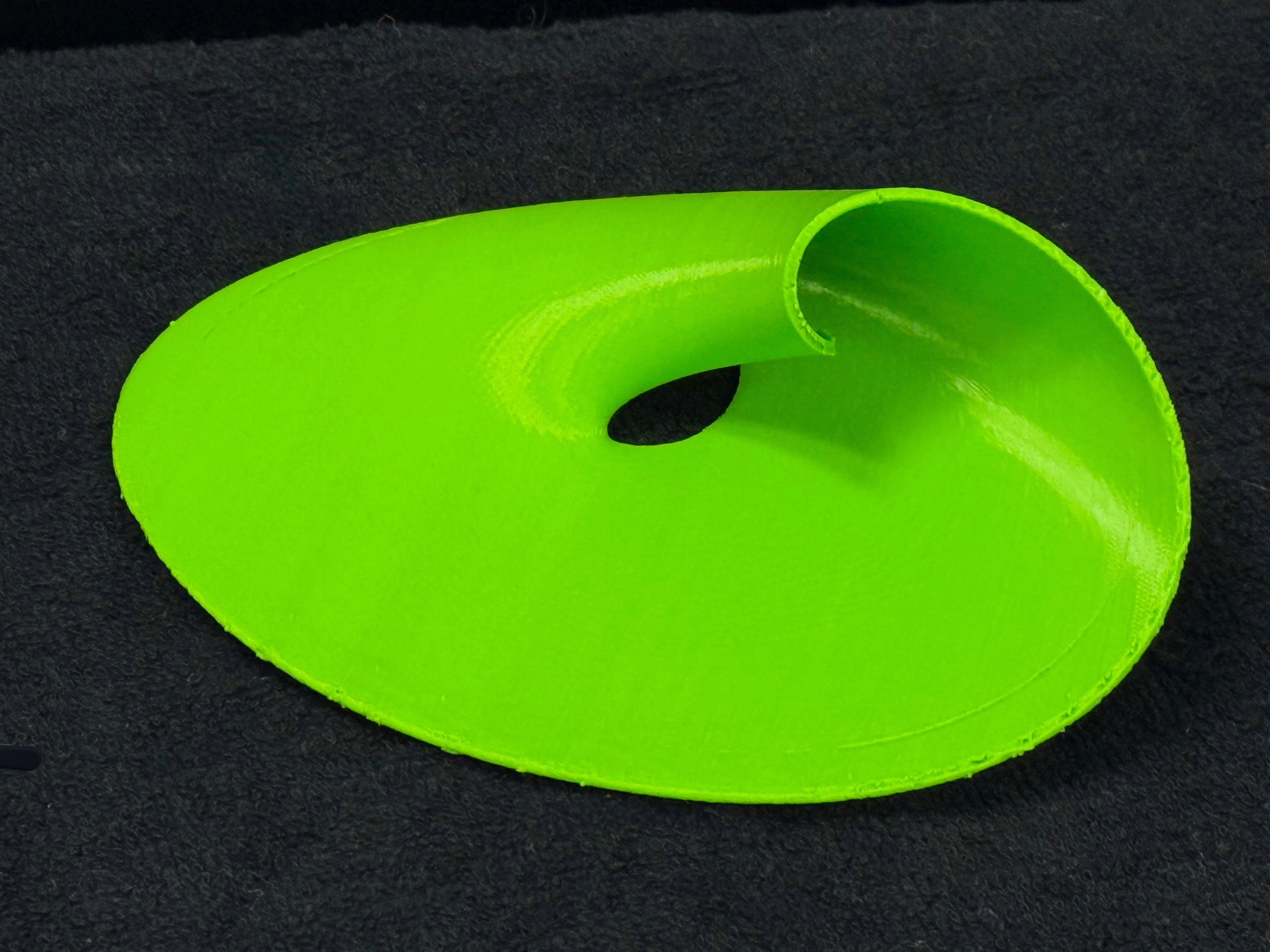}
\caption{Physical realization of the unstable manifold in the Arneodo-Coullet-Tresser system \eqref{eq:arneodo} with  $\beta = 0.4$ and $\mu = 0.863$. 
}
\label{fig:3dprint_manifold2}
\end{figure}

\begin{figure}[!htb]
\tabcolsep=0pt
\begin{center}
\begin{tabular}{cc}
\includegraphics[height=6cm, width=7cm, keepaspectratio]{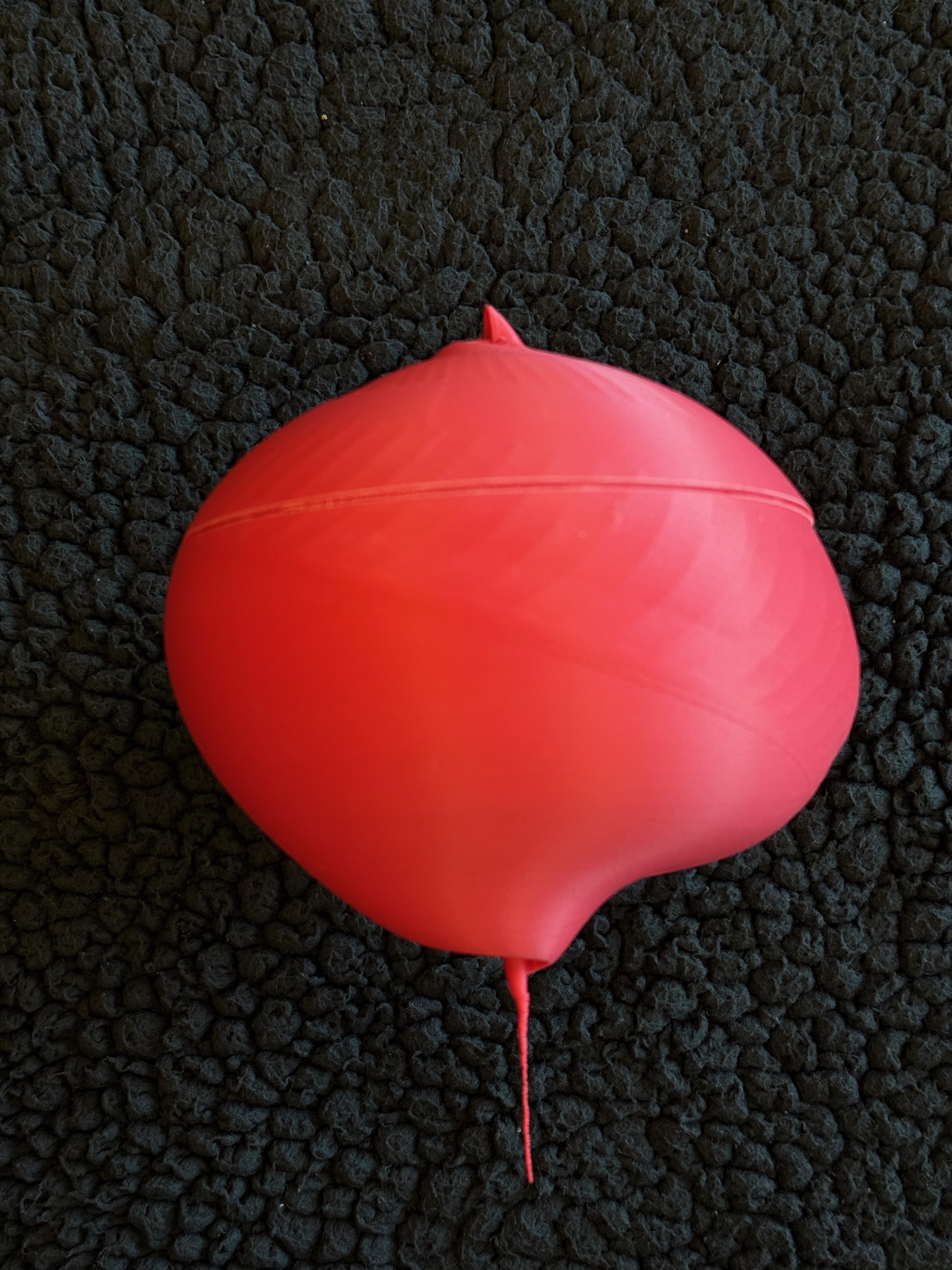} & 
\hspace{0.5em}\includegraphics[height=6cm, width=7cm, keepaspectratio]{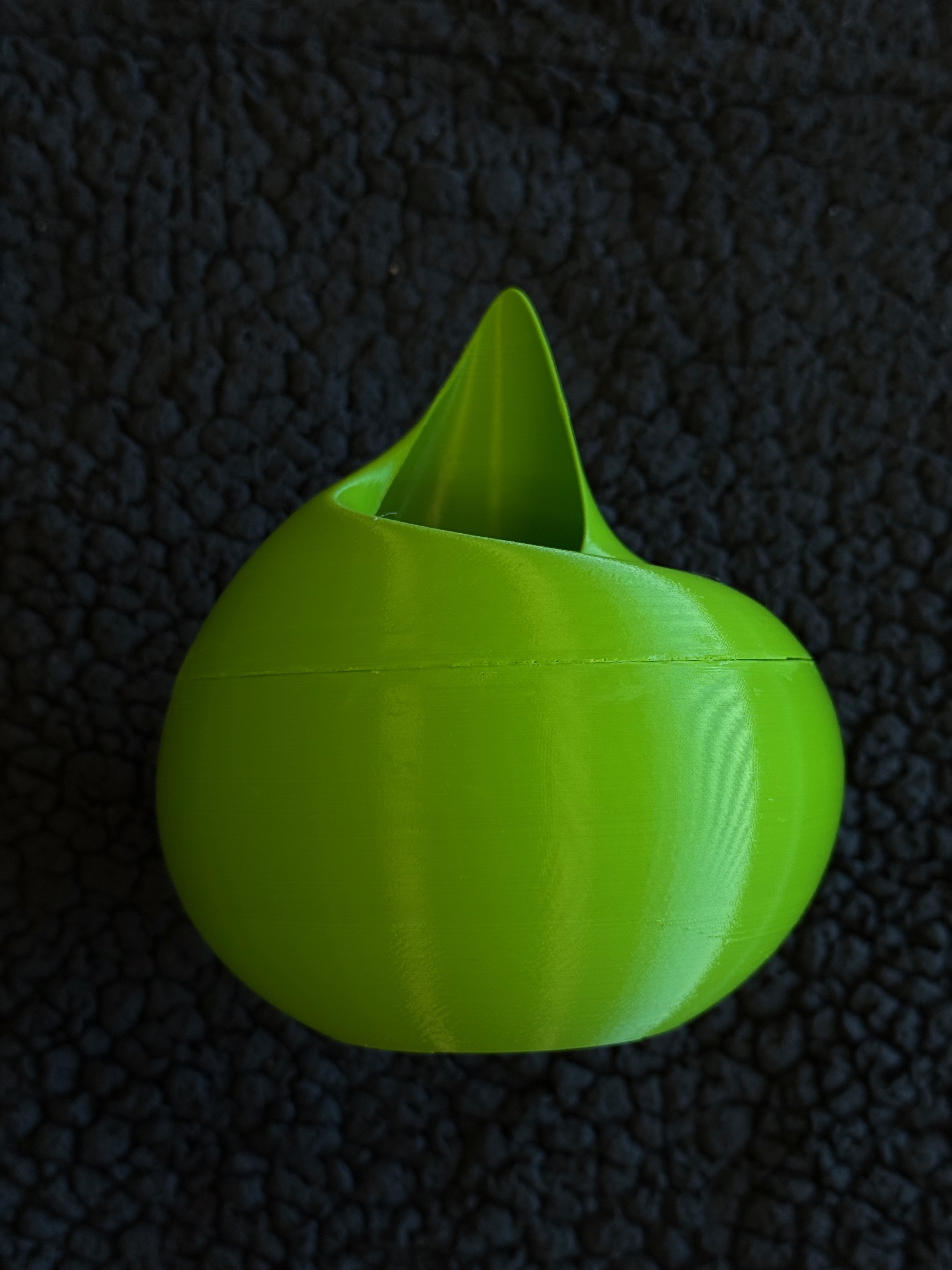} \\ 
\end{tabular}
\vspace{0.5cm} 
\begin{tabular}{c}
\includegraphics[height=6cm, width=7cm, keepaspectratio]{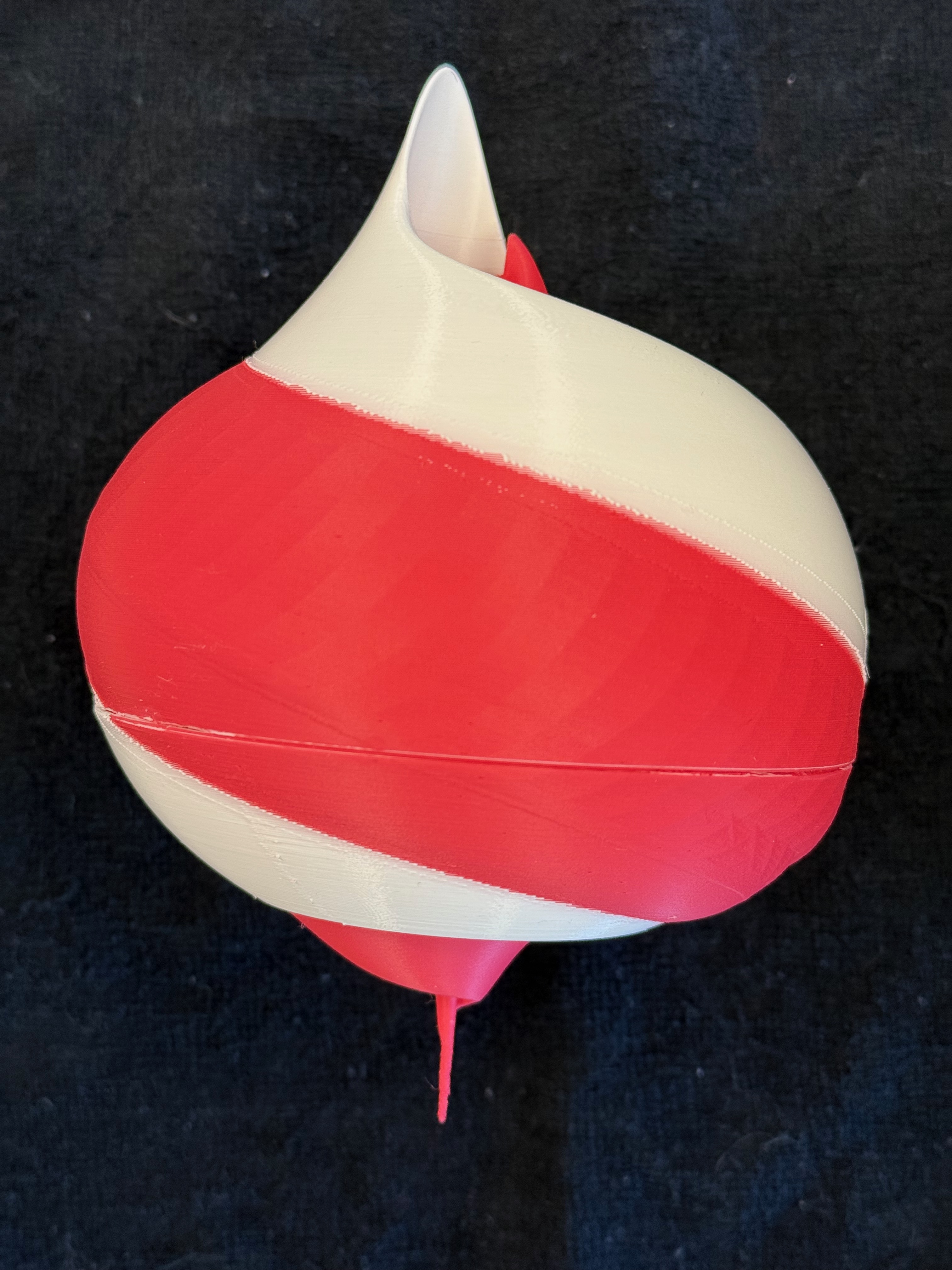} \\
\end{tabular}
\end{center}
\caption{Physical realization of the intersecting global manifolds in the Langford system. Top left: The 2D unstable manifold of the equilibrium point $p_1$. Top right: 2D stable manifold of the equilibrium point $p_2$. Bottom: A view highlighting the intersection between these two invariant manifolds, with a portion of one manifold protruding through the other. The parameter values for the Langford system used in these models are given in Section~\ref{sec:results}.}
\label{fig:3dprint_manifold3}
\end{figure}

This article introduces a novel approach to address this challenge: the application of 3D printing technology to create tangible, physical representations of geometric invariant sets, see Figures \ref{fig:phase_portrait_lorenz} --  \ref{fig:3dprint_manifold3}. We thus present a methodology for transforming abstract mathematical objects into accurate, three-dimensional models. This approach not only enhances our ability to visualize and study these intricate structures but also provides new pedagogical tools for teaching concepts in dynamical systems theory. Related work includes methods for 3D printing chaotic attractors~\cite{lucas20,bertacchini23,gagliardo20183d} and creation of crocheted invariant manifolds \cite{osinga2004crocheting}.

The process of transforming invariant manifolds for flows into physical 3D printed objects involves three distinct phases, each with its own computational and technical challenges. The first phase  focuses on generating an accurate numerical representation of the manifold. It consists of two parts: computing the local manifold, and computing the global manifold. For local manifold generation, we employ the Parameterization Method, a method which provides high-order approximations of the invariant manifold on a fundamental domain near the equilibrium solution. Computing the global manifold is done using integration schemes to extend the manifold to a desired size. If the eigenvalues happen to be complex conjugates of each other, this can be done by direct iteration, but otherwise this requires more sophisticated methods. We use the method of uniform arclength integration. 

The second phase addresses the crucial transition from a mathematical object to a printable 3D model. This requires implementing appropriate meshing schemes to create a triangulated surface representation. A key challenge here is the generation of a ``thickened" version of the manifold - converting the idealized surface of zero thickness into a physically printable object with appropriate thickness and structural integrity. This process must preserve the essential geometric features of the manifold while ensuring the resulting mesh is suitable for conversion to the file formats commonly used in 3D printing. The meshing algorithm must also account for potential self-intersections and ensure proper orientation of surface normals.

The final phase involves preparing the model for actual physical printing using slicing software. This stage requires careful consideration of multiple printing parameters: material selection, support structure generation for overhanging features, optimal build plate orientation to minimize supports while maintaining structural integrity, and layer height settings to balance print quality with production time. The choice of these parameters significantly impacts both the visual quality and the structural stability of the final printed manifold. Additionally, considerations must be made for the scale of the printed object, ensuring that fine geometric details are preserved while maintaining printability and handling requirements.

Our paper proceeds as follows. In order to lay the foundation for our methodology and results, we first review key concepts in the study of invariant manifolds in Section~\ref{sec:lin}. We proceed in Section~\ref{sec:comp} by looking at computational techniques for local manifolds using the Parameterization Method. Section \ref{sec:global} discusses our approach for the global manifold computation along with other tools used in the literature. Section~\ref{sec:3d} introduces our 3D printing pipeline, detailing the process of converting mathematical structures into printable models. Section~\ref{sec:results} demonstrates our approach through three examples: the stable manifold of the Lorenz system, the unstable manifold of the Arneodo-Cousset-Tresser flow, and the intersecting manifolds of the Langford system. We conclude in Section~\ref{sec:concl} with practical guidelines for reproducing our results,  providing references to our code repository in Section~\ref{sec:code}.

\section{Theory of Stable Manifolds.}\label{sec:lin}
A dynamical system describes evolution under a consistent set of rules. This very general 
class includes systems governed by principles in many different fields, including applications in physical, biological, chemical, and social sciences. In this paper, we restrict to the case 
in which the  state variable $x: \R \to \R^n$ is such that $x(t)$ is a vector-valued function of time, and our dynamical system can be written as a nonlinear ordinary differential equation of the form \[\dot{x} = f(x),\] where  $f: \R^n \to \R^n$ is a smooth function. 

In order to fix ideas, we start by considering planar linear systems \[\dot{x} = Ax,\]
where $x: \R \to \R^2$ is a vector-valued function, and $A$ is a $2$ by $2$ matrix. Note that if $x(0) = (0,0)$, then $\dot{x}(0) = (0,0)$, implying that $x(t)\equiv (0,0)$ for all $t$, which is an {\em equilibrium solution} for the system. 
Assume that $A$ has two eigenvalues $\lambda_1 < 0$ and $\lambda_2 > 0$ with corresponding eigenvectors $v_1$ and $v_2$. Then $Av_i = \lambda_iv_i$ for $i=1,2$, and all solutions are of the following form, where $c_1$ and $c_2$ are any real numbers:
\begin{align*}
x(t) = c_1e^{\lambda_1t}v_1 + c_2e^{\lambda_2t}v_2.
\end{align*}

\begin{figure}[!htb]
    \centering
    \includegraphics[width=1\textwidth]{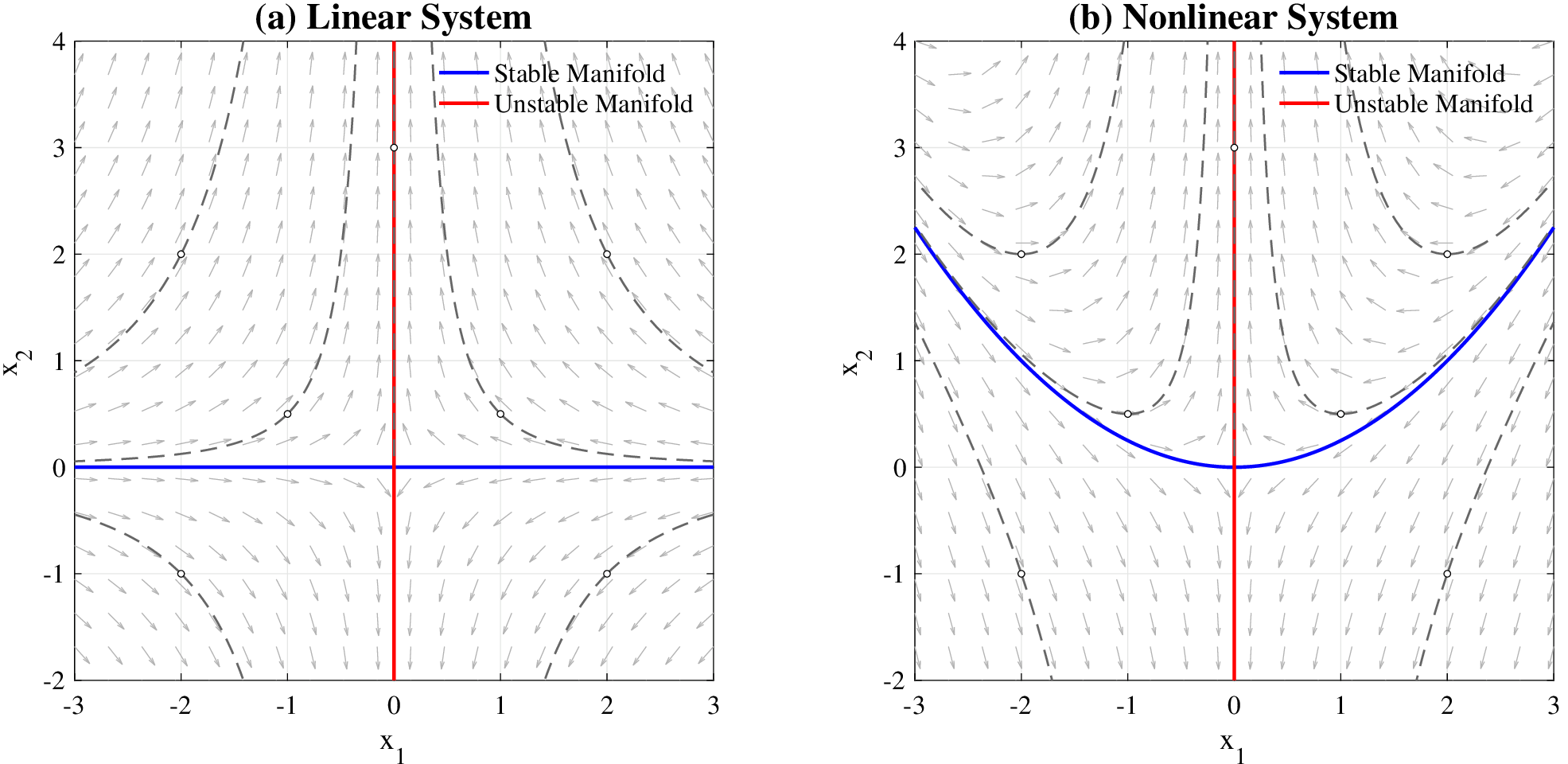}
    \caption{Phase portrait comparing linear and nonlinear dynamics: \textbf{(a)} Linear system $\dot{x}_1 = -x_1$, $\dot{x}_2 = 2x_2$ showing straight-line stable and unstable manifolds. \textbf{(b)} Nonlinear system $\dot{x}_1 = -x_1$, $\dot{x}_2 = 2x_2 - x_1^2$ with equilibrium $p = (0,0)$ and the same linearization as in (a), showing curved stable manifold of $p$ (blue), straight unstable manifold of $p$ (red), and sample trajectories (dashed). The light gray arrows indicate the vector field direction.}
    \label{fig:phase_portrait}
\end{figure}

These solutions are depicted in (Figure~\ref{fig:phase_portrait}a). Consider a particular solution $x_1(t)$ such that $x_1(0) = c_1 v_1 + c_2 v_2$ lies on the line $L_s$ through the origin in the direction $v_1$. This implies that $c_2 = 0$, and therefore $x_1(t)$ lies on $L_s$ for all values of $t$. Since any solution starting on $L_s$ stays on  $L_s$ for all time, we call $L_s$
an {\em invariant set} for the system. By a similar argument, the line $L_u$ that passed through the origin and is parallel to $v_2$ is also invariant. Note that the behavior on these two lines differs. In particular, since $\lambda_1<0$ on the line $L_s$, $\lim_{t \to \infty} x_1(t) =$ the origin. In fact, this uniquely characterizes $L_s$: there are no other points which limit to the origin. This line is the {\em stable manifold} of the origin.  In contrast, on $L_u$, since $\lambda_2>0$, $\lim_{t \to -\infty} x_2(t) =$ the origin, whereas for $t\to\infty$ the solution diverges exponentially quickly away from the origin. Again, there are no other points that limit to the origin as $t\to -\infty$. We call $L_u$ the {\em unstable manifold}. 

We have only considered a planar equation, but the same statement given above holds in higher dimensions: For the equation $\dot{x} = A x$, where $A$ is an arbitrarily large square matrix, suppose that $A$ has $n_s$ eigenvalues (with multiplicity) with real part less than zero, and $n_u$ eigenvalues (with multiplicity) with real part greater than zero. Then there is a 
$n_s$-dimensional stable linear subspace $E^s$ containing solutions that exponentially  decay to $0$ forward in time, and and 
$n_u$-dimensional unstable linear subspace $E^u$ with solutions that exponentially decay to zero backward in time. 

As a specific example of the situation above, when $A = \begin{pmatrix} -1 & 0 \\ 0 & 2 \end{pmatrix}$, $\lambda_1 = -1, v_1 = (1,0)$ and $\lambda_2 = 2, v_2 = (0,1)$, implying that  the stable and unstable manifolds are the $x$ and $y$ axes respectively (Figure~\ref{fig:phase_portrait}a). 

The fact that the stable and unstable manifolds are {\em linear subspaces} comes from the fact that the equation is linear. However, the  {\em stable manifold theorem} guarantees the existence of such sets for nonlinear equations.  That is, there is a local stable manifold, being a unique smooth invariant curve or surface of points converging to the origin exponentially quickly in forward time. There is also a similarly defined unstable manifold.  

Now consider the nonlinear differential equation
\begin{equation}\label{eq:non1}
\begin{split}
\dot{x}_1 &= -x_1 \\
\dot{x}_2 &= 2x_2 - x_1^2.    
\end{split}
\end{equation}
We can reasonably compare this to the linear example, since  $p= (0,0)$ is an equilibrium solution, and all the linear terms are the same as the linear example. 
A direct calculation shows that the general solution is given by  
\[x_1(t) = c_1 e^{-t}, \mbox{ and } x_2(t) =  \frac{c_1^2}{4}e^{-2t} + c_2e^{2t}, \]
for any $c_1,c_2 \in \R$. 
Define the two smooth curves
\[ W^s(p)  = \left\{(x_1,x_2): x_2 = \frac{x_1^2}{4}\right\} \mbox{ and } W^u(p) = \left\{(x_1,x_2): x_1 = 0 \right\}. \]
A direct check shows that both of these sets are invariant for the general solutions given above. Furthermore, all solutions in $W^s(p)$ limit to the origin as $t\to \infty$, and all solutions in $W^u(p)$ limit to the origin as $t \to -\infty$. Thus these are the stable and unstable manifolds of the origin (Figure~\ref{fig:phase_portrait}b). This is a satisfying result, but in general, it is not possible to find explicit solutions for a differential equation, and it is also
not possible to write down a simple formula for the stable and unstable manifolds. With this example in mind, in the next subsection, we give a formal treatments of invariant manifolds.

\subsection{Invariant Manifolds.} 

An invariant manifold is a subset of the state space of a dynamical system that remains invariant under the system's dynamics. In this entire section, we will use the following notation to describe our setup.

Consider the differential equation \begin{equation} \label{eq:ODE} \dot{x} = f(x),\end{equation} where $f:\R^n \to \R^n$ is a smooth function. 
Let $\phi_t(x)$ represent the {\em flow}: the solution of the differential equation at time $t$ with the initial condition $x$.
Assume $f(p) = 0$, implying that $x \equiv p$ is an equilibrium solution.
 The {\em linearization} of the differential equation at $p$, given by \begin{equation}\label{eq:linearized} \dot{y} = A y, \mbox{ where } \; A = Df(p).\end{equation} This linear equation has the closed form solution $e^{t A} y_0$, 
 where $y(0) = y_0$ is any initial condition. Here $e^{t A}$ is 
 the {\em matrix exponential}, defined for example in~\cite{meiss07}.
We refer to $p$ as a {\em hyperbolic} whenever $A$ has no eigenvalues with zero real part. 
%
We now define the stable and unstable manifolds of $p$.  
\begin{definition}[Stable and unstable manifolds]
 Assume that $p$ is a hyperbolic equilibrium for \eqref{eq:ODE}. 
 Define the following set of points in $\R^n$.
\begin{enumerate}
    \item Global stable manifold ($W^s(p)$): The set of all initial conditions $x$ such that $\phi_t(x)$ converges to $p$ for $t \to \infty$.
    \item Global unstable manifold ($W^u(p)$): The set of all initial conditions $x$ such that $\phi_t(x)$ converges to  $p$ for $t \to -\infty$.
\end{enumerate}
\end{definition}
The following theorem locally characterizes the stable and unstable manifolds\footnote{We assume a hyperbolic equilibrium, even though the theorem holds in a more general setting.}.
\begin{theorem}[Stable Manifold Theorem \cite{Guck}]
 Assume that $p$ is a hyperbolic equilibrium for \eqref{eq:ODE} with corresponding linearization \eqref{eq:linearized}. For matrix $A = Df(p)$, let $n_s$ and $n_u$ denote the number of eigenvalues (with multiplicity) with negative and positive real parts respectively.  Since $p$ is hyperbolic, $n_s + n_u = n$. Then for sufficiently small $\epsilon>0$, there exist unique invariant dynamical sets, known as {\em local} stable and unstable manifolds, defined as follows.  
{\small
\begin{equation}\label{localmanifolds}
    \begin{split}
W^s_{loc}(p) = &  \{ x \in W^s(p) :  \phi_t(x) \in B_\epsilon(p) \text{ for all } t \geq 0\}\\
 W^u_{loc}(p) = & \{ x \in W^u(p) :  \phi_t(x) \in B_\epsilon(p) \text{ for all } t \leq 0\},
    \end{split}
\end{equation}}
These local manifolds are surfaces through $p$ that are as smooth as $f$, with dimensions $n_s$, $n_u$ respectively. 
Let $E^s$ and $E^u$ be the stable and unstable subspaces corresponding to $\dot{y} = Ay$. At $p$, 
the local manifolds are tangent to $E^s$ and $E^u$ respectively. 
\end{theorem}
The above theorem uses the linearization \eqref{eq:linearized} to describe the local behavior of solutions converging to the equilibrium (for $t\to\infty$ or $t\to-\infty$).  In fact, we can describe the 
local behavior of {\em all} solutions starting near $p$, as stated in the following theorem. 
\begin{theorem}[Hartman-Grobman]
 Assume that $p$ is a hyperbolic equilibrium for \eqref{eq:ODE}, with corresponding 
 linearized equation \eqref{eq:linearized}. 
 Then there exists $\delta, s>0$ and a homeomorphism $h: B_s(0) \to B_\delta(p)$ such that  $h$ is a conjugacy between solutions of the linear equation and solutions of the nonlinear equation. 
 That is, for any $y \in B_s(0)$,
 \[
  \phi_t(h(y)) = h(e^{tA} y).
 \]
 This implies that $h$ maps trajectories of the linear equation to trajectories of the nonlinear equation, while preserving their topology and orientation. 
\end{theorem}
Now that we have a characterization of the local behavior, we comment on the global stable and unstable manifolds: They are obtained from the respective local manifolds by letting points in the local manifolds by direct integration. That is, 
\begin{align}\label{globalmanifolds}
    W^s(p) &= \bigcup_{t \leq 0} \phi_t(W_{loc}^s(p)) \\
    W^u(p) &= \bigcup_{t \geq 0} \phi_t(W_{loc}^u(p)). \nonumber
\end{align}
In the section, we give details on how to compute local manifolds. In the subsequent section, we explain how to use the local manifolds to compute the global manifolds. 

\section{Computation of Local Manifolds.}\label{sec:comp}

\subsection{The Parameterization Method.}\label{sec:param}
A powerful technique developed by Cabr\'{e}, Fontich, and de la Llave \cite{Cabre2003a} enables the computation of local stable and unstable  manifolds. In particular, the method involves finding a parameterization of the local manifold by using the invariance of the manifold. 

Before stating the general case, we give a basic understanding of the method by finding a parameterization of the local stable manifold by returning to the nonlinear example from Equation~\eqref{eq:non1} in Section~\ref{sec:lin}. 
Recall that \eqref{eq:non1} has a hyperbolic equilibrium $p=0$   and the linearized equation  has eigenvalues $\lambda_1 = -1 < 0$ and $\lambda_2 = 2 > 0$, with corresponding eigenvectors $v_1 = (1,0)$ and $v_2 = (0,1)$. Thus  $E^s = \mbox{ span}(v_1)$, the $x$-axis, and $E^u$ is the $y$-axis. By the stable manifold theorem, we know that $W^s_{loc}(p)$ is a one-dimensional curve that intersects $p=0$, and is tangent to the $x$-axis at $0$. 
By the Hartman-Grobman Theorem, we know that there exists $h:B_s(0)\to B_\delta(p)$ such that  $h(x\mbox{-axis})= W^s_{\rm loc}(p)$. We use this to define a parameterization $P(\theta)$ for the stable manifold as follows.  
For every $|\theta|<1$, define \[P(\theta) = h((\theta,0)) .\] 
We now use the Hartman-Grobman Theorem to find an equation for $P$. Since $\lambda_1$ is $-1$ we know that the linear solution starting at $(\theta,0)$ has value $e^{t A} = (e^{-t} \theta,0)$ at time $t$.
Therefore the conjugacy equation for $h$ implies
\[
\phi_t(P(\theta)) = P(e^{-t}\theta).
\]
Differentiating the equation  with respect to $t$, setting $t=0$, and rearranging the terms, we get 
\begin{equation}\label{conjugacy}
f(P(\theta)) =  -\theta DP(\theta).
\end{equation}
We now write $P(\theta)$ in the form of a Taylor series 
\[P(\theta) = \begin{pmatrix} a_{0} + a_{1}\theta + a_{2}\theta^2 + \cdots \\ b_{0} + b_{1}\theta + b_{2}\theta^2 + \cdots \end{pmatrix}\]
where $a_{0} = b_{0}= 0$ since $p=0$ is contained in the stable and unstable manifolds. 
We find the rest of the coefficients by substituting $P$ into equation~\ref{conjugacy}, giving
\begin{align*}
-\theta DP(\theta) &= \begin{pmatrix} -\theta(a_1 + 2a_2\theta + 3a_3\theta^2 + \cdots) \\ -\theta(b_1 + 2b_2\theta + \cdots) \end{pmatrix} \\
f(P(\theta)) &= \begin{pmatrix} -(a_1\theta + a_2\theta^2 + \cdots) \\ 2(b_1\theta + b_2\theta^2 + \cdots) - (a_1\theta + a_2\theta^2 + \cdots)^2 \end{pmatrix}
\end{align*}
Solving term-by-term, we find that:
\begin{itemize}
\item For the first component: $a_1$ is arbitrary (we choose $a_1 = 1$), and $a_k = 0$ for $k > 1$.
\item For the second component: $b_1 = 0$, $b_2 = \frac{1}{4}$, and $b_k = 0$ for $k > 2$.
\end{itemize}

Therefore, $P_1(\theta) = \theta$ and $P_2(\theta) = \frac{1}{4}\theta^2$, revealing that the stable manifold is given by:
\begin{equation}
x_2 = \frac{1}{4}x_1^2,
\end{equation}
as stated when we first introduced the equation. 

Returning to the general theory, consider \eqref{eq:ODE} with a hyperbolic equilibrium point $p$, where $f: \mathbb{R}^n \to \mathbb{R}^n$ is a real analytic vector field and  $Df(p)$ is diagonalizable\footnote{The diagonalizability assumption simplifies the recursive computation of coefficients but is not essential as Newton method may be used.} 
with eigenvalues $\{\lambda_i^s\}_{i=1}^{n_s}$ such that for all $i$, $Re(\lambda_i^s)<0$ and  $\{\lambda_i^u\}_{i=1}^{n_u}$ such that for all $i$, $Re(\lambda_i^u)>0$, such that the eigenvalues are distinct (and satisfy certain technical nonresonance conditions). Let $\{\xi_i^s\}_{i=1}^{n_s}$, $\{\xi_i^u\}_{i=1}^{n_u}$ be the corresponding eigenvectors.
As before, let $h$ be the conjugacy function given by the Hartman-Grobman Theorem, and let $\theta \in \R^{n_s}$ with $|\theta|\le 1$. Consider the vector $v_s$ in  the stable linear subspace  $E^s$ given by \[v_s = \theta_1 \xi^1_s + \dots + \theta_{n_s} \xi^{n_s}_s.\] Then using the scale parameter $s$, define the parameterization $P$ of the stable manifold by 
\[
P(\theta) = h(s v_s).
\]
This gives a parameterization $P: B_1(0) \subset \mathbb{R}^{n_s} \to \mathbb{R}^n$ for $W_{\textrm{loc}}^s(p)$.
Using the same argument as before, differentiating this equation with respect to $t$ and setting $t=0$, we see that 
\begin{equation}\label{eq:3}f(P(\theta)) = DP(\theta)\Lambda \theta,\end{equation} 
where  $\Lambda$ is the diagonal matrix with the eigenvalues $\lambda^i_s$, $i=1,\dots, n_s$ on the diagonal. In order to approximate $P$, we write $P(\theta)$ as a series solution, and use equation~\eqref{eq:3} to solve for the coefficients in the series. In the previous example, we were able to compute an exact formula for the stable manifold, but in general 
we can only approximate the series for $P$ with a finite number of terms. 
In particular we refer to the truncated series $P(\theta) = \sum_{j,k = 1}^N P_{jk} \theta_1^j \theta_2^k$ as an {\em order $N$ approximation.}
For a more detailed explanation and implementation of the Parameterization Method, one can refer to \cite{fleurantin2020resonant}. 

In order to be able to 3D print our results,  we are focusing on the case of a three-dimensional phase space. In this case, the only possibilities are that the equilibrium 
has stable (resp. unstable) manifolds of dimension $0(3),1(2),2(1)$, or $3(0)$. The first and last case are uninteresting, as one of the two manifolds is locally the whole space. Therefore we only consider the cases with the stable and unstable manifolds are dimensions 1 and 2 or 2 and 1. Computing the  one-dimensional manifold is completely analogous but more  straightforward to calculate than the two-dimensional one, and therefore we give a detailed description on computing the two-dimensional manifold.\footnote{Code is provided for both the one- and two-dimensional manifolds.} 
A general algorithm\footnote{The notation used Algorithm \ref{alg:loc1} (indices, $r$, $\theta$) aligns with the implementation provided in the accompanying GitHub repository to facilitate direct comparison.} to compute the local stable manifold is outlined in Algorithm \ref{alg:loc1}. The unstable manifold computation is completely analogous; note that the two manifolds must be computed separately. In Steps 13--15 of this Algorithm, we evaluate $P(\theta)$ at the points furthest from the equilibrium $p$, i.e. on the circle $C_0 = \{\theta: |\theta-p| = 1 \}$. This topological circle is the set of points that we use 
in order to compute the global stable manifold in the next section. 
\begin{algorithm}
\caption{Parameterization Method for two-dimensional stable manifolds, order $N$}
\begin{algorithmic}[1]
\Require System parameters
\Ensure Local stable  manifold parameterization
\State Compute hyperbolic equilibrium point $p$
\State Define the righthand side of the differential equation using $f(\zeta)=Df(p)\zeta + R(\zeta)$
\State Find stable  eigenvalues $\lambda_1$, $\lambda_2$ with $Re(\lambda_1) \le Re(\lambda_2)$. Find corresponding eigenvectors $v_1$, $v_2 \in \mathbb{C}^3$
\State Initialize the coefficients $P_{j,k}$ of order $N$
\State Set $P_{0,0} = p$, $P_{1,0} = sv_1$, $P_{0,1} = sv_2$ \Comment{$s$ to bound order $N$ coefficients}

\For{$k = 2$ to $N$}
    \For{$j = 0$ to $k$}
        \State $n_1 \gets k-j$, $n_2 \gets j$
        \State Solve $Df(p)P_{n_1,n_2} - (n_1\lambda_1 + n_2\lambda_2)P_{n_1,n_2} = R_{n_1,n_2}$ 
    \EndFor
\EndFor
\Ensure Evaluate the  parameterization approximation on the unit circle 
\For{each point $(\sigma_1, \sigma_2) = (r\cos\theta_i, r\sin\theta_i)$}
    \State Compute $P(\sigma_1, \sigma_2) = \sum_{m_1 + m_2 \leq N} P_{m_1,m_2}(\sigma_1 + i\sigma_2)_1^{m_1}(\sigma_1 - i\sigma_2)^{m_2}$
\EndFor

\Return Manifold points $P(\sigma_1, \sigma_2)$

\end{algorithmic}\label{alg:loc1}
\end{algorithm}

\begin{remark}
The algorithm assumes sufficient smoothness of the vector field. The accuracy of the approximation improves as the order $N$ increases, subject to the radius of convergence of the series.
\end{remark}

\begin{remark}
The choice of scale parameter $s$ affects the size of the computed local manifold patch. Optimal selection of $s$ in line 5 of Algorithm~\ref{alg:loc1} depends on the specific system properties and desired visualization region.
\end{remark}

\section{Computation of Global Manifolds.}\label{sec:global}
Once a  local manifold has been determined, we need to integrate to obtain the global manifold. 
Many different numerical methods for computing global manifolds are described in the comprehensive book of Krauskopf, Osinga, and Galán-Vioque~\cite{Krauskopf2007}, papers of Dellnitz and Junge~\cite{dellnitz2016computation}, Krauskopf et al.~\cite{Krauskopf05}, Henderson~\cite{henderson2005computing}, and Haro~\cite{haro16}.

The general method for finding the global manifolds is to use the definition of the global manifolds in terms of the local manifolds given in Equation \eqref{globalmanifolds}. Note that the global manifolds are in general unbounded, so the best we can do is to integrate these equations up to some distance from the equilibrium.  
In some cases, this simply involves starting on the boundary of the local stable (resp. unstable) manifold, and integrating backward in time for $t$ from $0$ to $-T$ to get an additional portion of the global stable manifold (resp. forward in time for $t$ from $0$ to $T$ for an additional portion of the global unstable manifold). 

Specifically, let $C_0$ be the circle of points furthest from $p$ on the local manifold. In order to ensure that the points are well spaced, we reparameterize $C_0$ with respect to the arclength parameterization given by the curve $r_0(\eta)$ where  $\eta$ varies from $0$ to $1$ and $r_0(0) = r_0(1)$. We wish to integrate 
this curve backwards (forwards) in time. That is, we wish to solve the following ordinary differential equation, where we allow the solution $y$ to depend on the independent variable $t$ but also on the variable $\eta$, as this allows us to parameterize the set of initial conditions. 
\[
\frac{d}{dt} (y(t,\eta)) = f(y(t,\eta)) \; , \mbox{ where } y(0,\eta) = r_0(\eta). 
\]
While it is possible to write a numerical solver to solve this, since our goal is to create a printable mesh, we can also take advantage of a built-in solvers that represents solutions as parameterized surfaces, meaning that we can take advantage of built-in meshing algorithms. 

\begin{figure}[!htb]
    \centering
    \begin{minipage}[b]{0.55\textwidth}
        \centering
        \textbf{(a)} \\
        \includegraphics[width=\textwidth]{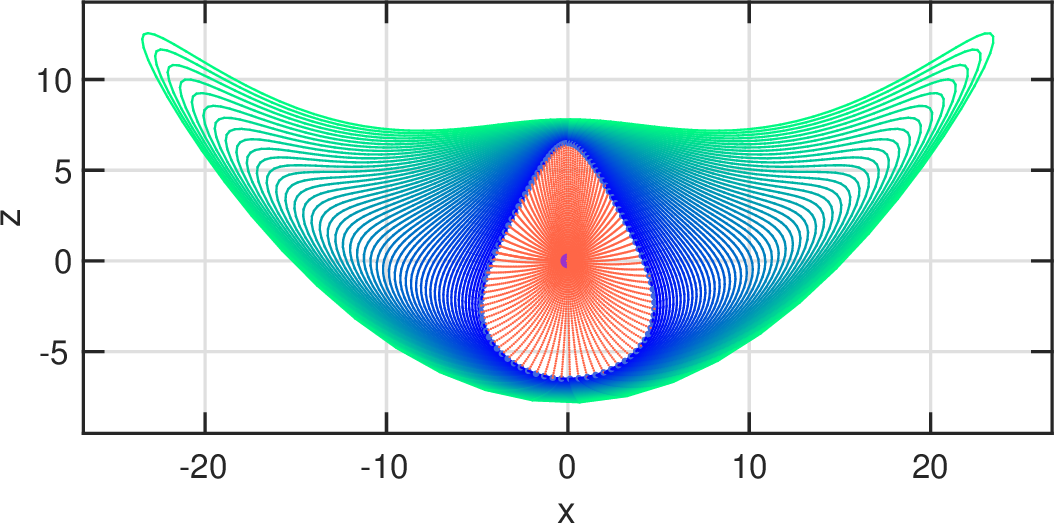}
        \label{fig:stretching}
    \end{minipage}
    \hfill
    \begin{minipage}[b]{0.55\textwidth}
        \centering
        \textbf{(b)} \\
        \includegraphics[width=\textwidth]{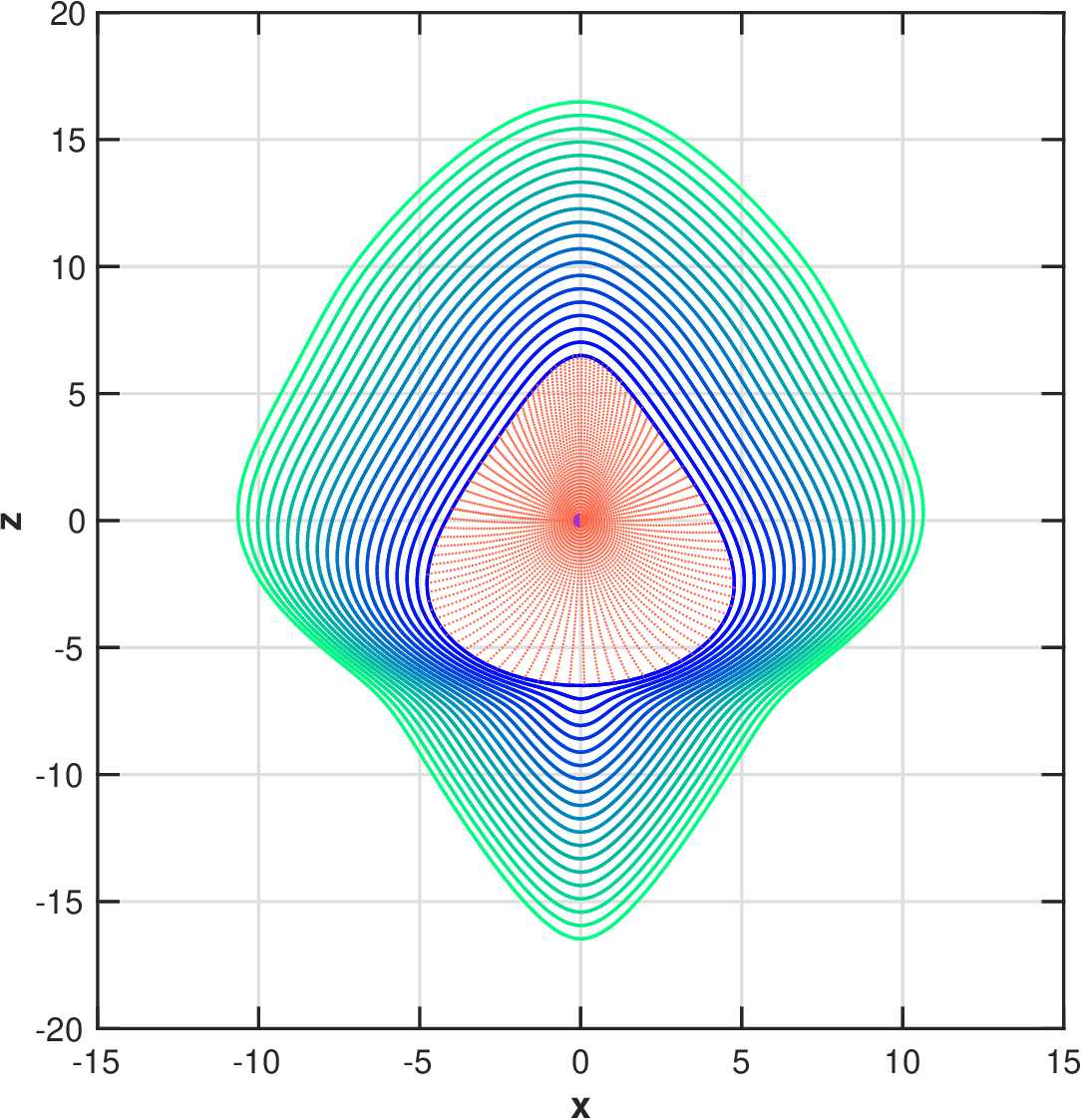}
        \label{fig:arclength}
    \end{minipage}
    \caption{Comparison of integration methods for computing the global stable manifold of the Lorenz system~\eqref{eq:lorenz} In both figures, the red inner region in the center represents the local stable manifold $W^s_{loc}(p)$ directly surrounding the equilibrium point at the origin.
    \textbf{(a)} When the eigenvalues $\lambda_1 < \lambda_2$ are real and distinct, direct integration leads to unequal stretching. \textbf{(b)} The arclength parameterization method in Algorithm~\ref{alg:p2} solves this problem.}
    \label{fig:manifold_comparison}
\end{figure}

If the eigenvalues $\lambda_1$ and $\lambda_2$ are complex conjugates, then they have the same stretching factor, and the method of direct integration works to give a nice depiction of the global manifold without further modification.  However, if $\lambda_1 < \lambda_2$, this method lead to unequal stretching during the integration process, and as $T$ grows the stretching becomes so extreme that our integration essentially only captures a single dimension of the manifold, see Figure \ref{fig:manifold_comparison}a.

A solution to this stretching problem~\cite{Krauskopf05} is to integrate in such a way that the time is parameterized by arclength, see Figure \ref{fig:manifold_comparison}b. Assume the original differential equation is given by  \[ \dot{x} = f(x), \mbox{ where } x = (x_1,x_2,x_3), \; f(x) = (f_1(x),f_2(x),f_3(x)), \] denote  \[ \|f\|_2 = \sqrt{(f_1(x))^2+(f_2(x))^2 + (f_3(x))^2}\; .\] Our reparameterized differential equation   is given by  
\begin{equation}\label{eq:arclength}
\frac{dx}{d\tau} = \frac{f(x)}{\|f\|_2}, \; \;\;
\frac{dt}{d \tau} = \frac{1}{\|f\|_2}\; .
\end{equation}

This yields the same trajectories as the original equation, but the time rescaling implies that along the curve in the global manifold with a fixed value of $\tau$, all the points are the same geodesic distance from the starting curve $r_0$. In addition to this modification of our original vector field, if we calculate with the initial conditions for too long, the 
resulting points will start to bunch up in small regions, thereby not representing the full behavior along a topological circle. Therefore, we break our full length $T$ into $K$ 
time steps of length $T/K$. Before any time step from $T_j$ to $T_{j+1}$, we find a reparametrization of the outermost points with even spacing  via arclength, i.e. of the computed circle of points at time $T_{j}$. This method of calculation is summarized in Algorithm \ref{alg:p2}.

\begin{algorithm}[tbh!]
\caption{Two-dimensional global stable/unstable manifold computation} 
\begin{algorithmic}[1]
\Require 
\State Initial data points $\mathcal{D}$ from local manifold computation
\State System parameters $\Theta = \{\theta_1, \ldots, \theta_k\}$
\State Time interval $[0,T]$
\State Parameter space $\Omega = \{(\sigma_1,\sigma_2) : \|\sigma\| \leq r\}$ \Comment{2D disc}
\State Numerical precision parameters $\epsilon_{tol}$
\Ensure Two-dimensional global invariant manifold $\mathcal{M}$ and its parameterization

\Procedure{2DManifoldComputation}{}
    \State \textbf{Phase 1:} System Initialization
    \State Initialize boundary curve from $\mathcal{D}$, evenly spaced in $\R^3$
    \State Define vector field $F: \mathbb{R}^n \to \mathbb{R}^n$, \Comment{For $n=3$ or $n=4$ (in the case of \eqref{eq:arclength})}
    
    \State \textbf{Phase 2:} Differential System Setup
    \For{boundary $\partial\Omega$}
        \State Define $\dot{X} = F(X)$ 
        \State Set IC: $X(0)$ \Comment{Set either pointwise or as an equation}
    \EndFor
    
    \State \textbf{Phase 3:} Numerical Integration
    \State Choose integration method $\mathcal{I}$ (e.g., NDSolve, Method of lines )
    \For{each time step $T_j \in [0,T]$}
        \For{along boundary $ \partial\Omega$}
            \State Update solution: $Y(t,x), T_j<t \le T_{j+1}$
            \State Reparameterize $Y$ at $T_{j+1}$ to be an evenly spaced topological circle
        \EndFor
        \State Check convergence: $\|X(t_j) - X(t_{j-1})\| < \epsilon_{tol}$ \Comment{may be built in methods with user set tolerances}
    \EndFor
    
   \EndProcedure

\Return Two-dimensional manifold $\mathcal{M}$ suitable for 3D printing
\end{algorithmic}\label{alg:p2}
\end{algorithm}

\section{3D printing invariant manifolds.}\label{sec:3d}

Having established the basic principles and general computational methods for generating invariant manifolds, we now turn to the practical challenges of transforming these mathematical structures into physical objects.

Recent advances in additive manufacturing technology have revolutionized the creation of mathematical objects~\cite{shahrubudin2019overview,lee2017fundamentals}. From industrial prototyping to educational visualization~\cite{segerman20123d}, new techniques and materials continue to enable increasingly sophisticated prints~\cite{zastrow2020new}. These developments, coupled with the standardization of 3D printing workflows~\cite{jandyal20223d}, provide the foundation for tackling the unique challenges of printing invariant manifolds.

However, despite these technological advances, the specific task of converting computed manifold data into 3D-printable models requires careful consideration of both numerical accuracy and the physical constraints of additive manufacturing. Our methodology addresses several key challenges, including mesh/surface generation, structural stability, and the preservation of mathematical features through the printing process. In what follows, we detail our complete pipeline from numerical data to finished print, focusing primarily on two-dimensional invariant manifolds of three-dimensional vector fields, with particular attention to the technical requirements that ensure both mathematical fidelity and practical printability.

The main problem that we need to overcome is that a parameterized surface is not actually a solid, in that it contains no volume. In order to be able to print a surface, we need some way to thicken it so that it can print. There are various ways to do this, but we have found that the most reliable way to do it is by using a Matlab code which performs the Parameterization Method; followed by bringing the parameterization into  Mathematica in order to generate the object using the built in integration methods and built in thickness commands for parameterized surfaces. In addition, Mathematica allows us to create thickened tubes around curves within the invariant manifolds, in order to be able to see the internal structures that would otherwise be hidden by the surface. In addition, we have also had some luck in removing  a ``pie slice" cut into the object to show the internal structure. See Figure~\ref{fig:phase_portrait_langford}.

There are many other ways that one might consider performing the calculations, and we have in fact tried several other methods, with quite unsatisfying results. For example, we tried using built in thickening methods in the software Blender to create the solid object, but invariant manifolds are sufficiently complex objects that in general we end up with a mesh file with significant structural errors that make it unprintable. We also tried to do the entire calculation in Matlab, computing the object using triangulation mesh and solidification methods that have been created for Matlab. Again here, the problem is that although in theory it appears to create a perfectly meshed object, we often find that it does not successfully print due to holes and incorrect normals as a result of meshing issues. 

\section{3D Printed Manifold Examples}\label{sec:results}
To demonstrate the effectiveness of our computational approach and visualization pipeline, we examine three dynamical systems of particular interest. These systems serve as compelling test cases, illustrating both the robustness of our method and its ability to reveal intricate geometric structures through three-dimensional visualization and physical reproduction.

\subsection{Stable Manifold of the Lorenz System.}
The Lorenz system is a three-dimensional system of ordinary differential equations defined by the vector field
\begin{equation}\label{eq:lorenz}
    \begin{aligned}
        \dot{x} &= \sigma(y-x), \\
        \dot{y} &= \varrho x - y - xz, \\
        \dot{z} &= xy - \beta z,
    \end{aligned}
\end{equation}
with standard parameters $\sigma = 10$, $\varrho = 28$, and $\beta = 8/3$, where we consider the equilibrium at the origin $p=(0,0,0)$. Originally developed by Lorenz (1963) as a simplified model of atmospheric convection, this system is a classical example of sensitive dependence on initial conditions, where nearby trajectories diverge dramatically under the flow. For the standard parameter values, trajectories converge to the famous butterfly attractor, also known as the Lorenz attractor, which occupies a small region in $\mathbb{R}^3$. This attractor is particularly intriguing from a geometric perspective, as it exists between a one-dimensional and two-dimensional object, possessing a fractal dimension that characterizes it as a strange attractor, making it a fundamental example in the study of chaotic dynamical systems. 

To visualize the stable manifold of the Lorenz system at $p$ with classical parameters $\sigma = 10$, $\varrho = 28$, and $\beta = 8/3$, we first linearize, 
noticing that there are two stable eigenvalues $\lambda_1 \approx -22.8 < \lambda_2 \approx -2.7$  (as well as one unstable eigenvalue $\lambda_3 \approx 11.8 >0$). The fact that the two stable eigenvalues are real and different in size implyes that we require an arclength parameterization method for our global stable manifold computation. We 
employ Algorithm~\ref{alg:loc1} to compute the local stable manifold through iterative solution. 

The resulting coefficients of the parameterization enable us to generate a discretized boundary curve. We then use Algorithm~\ref{alg:p2} using the arclength parameterization of the vector field $F$ given in \eqref{eq:arclength}. Figure~\ref{fig:manifold_comparison}a shows how stretched the manifold would be if we did not use the arclength parameterization. The use of Algorithm~\ref{alg:p2} allows us to compute the global manifold structure up to $T = 180$ with $100$ total intermediate time steps. The result 3D print is shown in Figure~\ref{fig:3dprint_manifold}.  

\subsection{Unstable Manifold of the Arneodo-Coullet-Tresser System.}
The Arneodo-Coullet-Tresser system describes a three-dimensional dynamical system that exhibits Shilnikov chaos. The system is defined by the vector field
\begin{equation}
\begin{aligned}
\dot{x} &= y, \\
\dot{y} &= z, \\
\dot{z} &= -y - \beta z + \mu x(1-x),
\end{aligned}\label{eq:arneodo}
\end{equation}
with the parameters $\beta = 0.4$ and $\mu = 0.863$ where Shilnikov chaos is known to occur. Originally proposed as a simple model demonstrating spiral chaos, this system is particularly notable for its explicit exhibition of Shilnikov-type homoclinic orbits. The system possesses two equilibria: a saddle-focus at the origin $p_1=(0,0,0)$ and a second equilibrium point $p_2=(1,0,0)$ that undergoes various bifurcations as parameters change, ultimately leading to the formation of a Shilnikov attractor.

To visualize the Shilnikov attractor of the Arneodo system at $\beta = 0.4$ and $\mu = 0.863$, we note that $DF(p_2)$ has unstable eignevalues $\lambda_{1,2} \approx 0.1542 \pm 1.0930i$. Since these are complex conjugates, we will not need to use the arclength parameterization for our global manifold calculations. 
We first employ Algorithm~\ref{alg:loc1} to compute the local manifold structure 
The resulting parameterization coefficients enable us to generate a discretized boundary curve, which Algorithm~\ref{alg:p2} then evolves into the global manifold structure up to $T = 22$ using NDSolve which typically employs adaptive time stepping. The resulting 3D printed manifold is shown in Figure~\ref{fig:3dprint_manifold2}.

\subsection{Stable and Unstable Manifolds of the Langford System.} The Langford system describes a three-dimensional dynamical system that exhibits both Hopf and cusp bifurcations. The system is defined by the vector field
\begin{equation}
\begin{aligned}
\dot{x} = &
(z-\beta)x - \delta y,\\
\dot{y} = & \delta x + (z-\beta)y, \\
\dot{z} = & \tau + \alpha z - \frac{z^3}{3} - (x^2 + y^2)(1+\varepsilon z) + \zeta zx^3,
\end{aligned}\label{eq:langford}
\end{equation}
with classical parameters $\varepsilon = 0.25$, $\tau = 0.6$, $\delta = 3.5$, $\beta = 0.7$, $\zeta = 0.1$, and bifurcation parameter $\alpha > 0$. Originally derived by truncating a normal form of a simultaneous Hopf/cusp bifurcation to second order, the system includes an additional third-order term that breaks axial symmetry. The system serves as a model for dissipative vortex dynamics and rotating viscous fluids, featuring a $z$-axis invariant subsystem and rich dynamics including periodic orbits and multiple equilibrium solutions.

To visualize the intersection of global invariant manifolds in the Langford system
with parameters $\alpha = 1.1022$, $\beta = 0.7$, $\delta = 3.5$, $\gamma = 0.6$, $\zeta = 0.1$, and $\varepsilon = 0.25$, we
calculate that the system has three equilibria on the $z$-axis: $p_1 \approx (0,0,2.05)$, $p_2 \approx (0,0,-1.43)$, and $p_3 \approx (0,0,-0.61)$. For the equilibrium $p_1$, the linearization $DF(p_1)$ has eigenvalues 
$\lambda_{1,2} \approx 1.35 \pm 3.50i$ and $\lambda_3 \approx -3.08$, corresponding to a saddle-focus with a two-dimensional unstable manifold and a one-dimensional stable manifold. For the equilibrium $p_2$, the linearization $DF(p_2)$ has eigenvalues 
$\lambda_{1,2} \approx -2.13 \pm 3.50i$ and $\lambda_3 \approx -0.95$, indicating a sink with complex conjugate eigenvalues. The presence of complex conjugate eigenvalues implies that we will not need to parameterize with respect to arclength when computing the
global manifold structure. 

Now we employ Algorithm~\ref{alg:loc1} to compute the local manifolds of two distinct equilibrium solutions on the $z$-axis. We then use the resulting discretized boundary curves and evolve using time integration to capture the global manifold structures using Algorithm~\ref{alg:p2}. 

Figure~\ref{fig:3dprint_manifold3} shows intersections of the resulting structures: the global unstable manifold of $p_2$ with $T = 9.5$ and the global stable manifold of $p_3$  $T = 2.3$. 
We see that these two manifolds intersect  transversally  along curves. The integration was performed using Mathematica's NDSolve with adaptive time-stepping to ensure numerical accuracy. 

Figure~\ref{fig:phase_portrait_langford} shows 3D printed unstable manifolds for two different $\alpha$ parameter values. The yellow filamentary structure represents the 2D unstable manifold of the equilibrium point $p_1 \approx 1.84$ for $\alpha=0.806$, while the purple surface shows the 2D unstable manifold of $p_1 \approx 1.94$ for $\alpha=0.95$. The yellow invariant manifold, shown from two different perspectives in the right images, was computed using a numerical integration method that approximates the global invariant structure by iteratively applying the flow to an initial set of points along the local manifold, with interpolation between successive iterations. For the purple manifold (left image), set difference techniques were employed to remove extraneous portions, allowing better visualization of its complex geometric structure.


\section{Conclusions and how to make your own.}\label{sec:concl}
This article has presented a comprehensive computational pipeline for visualizing and physically realizing invariant manifolds of dynamical systems through 3D printing. Our approach combines careful mathematical computation with modern manufacturing techniques, enabling both theoretical study and tactile exploration of these complex geometric structures. The methodology has been successfully demonstrated by computing and visualizing invariant manifolds associated with three dynamical systems: the stable manifold at the origin of the Lorenz system, the unstable manifold of a saddle-focus for the Arneodo-Coullet-Tresser system, and the intersecting stable and unstable manifolds of equilibrium solutions in the Langford system, showcasing the versatility and robustness of our approach.

\subsection{Implementation Framework.} The reproduction of these results begins with careful system definition and analysis. For any dynamical system of interest, one must first clearly define the vector field and identify relevant equilibrium solutions. The choice of system parameters significantly influences the manifold structure and should be selected based on the specific phenomena one wishes to study. Our computational pipeline then proceeds with the local manifold computation using Algorithm~\ref{alg:loc1}, implemented in our provided Wolfram Mathematica\textsuperscript{\textregistered} scripts. This crucial step requires attention to numerical parameters, particularly the order of parameterization $N$ and the scaling factor, which should be adjusted based on the system's properties and desired accuracy. 

Global manifold generation follows through Algorithm~\ref{alg:p2}, where the local manifold boundaries are evolved to capture the complete invariant structure. This phase demands careful monitoring of numerical stability during time integration and may require adaptive step-size control. The resulting geometric data undergoes appropriate transformations to ensure optimal visualization while maintaining mathematical accuracy. The process concludes with the preparation for 3D printing, where computed surfaces are exported to standard formats (.stl or .obj) and processed through slicing software such as Cura to generate the necessary G-code for physical realization. That is, the slicer creates precise instructions to the printer for how hot to heat the filament, the path of motion, and where to extrude filament.

\subsection{Practical Considerations.} As discussed in Section \ref{sec:3d}, several challenges commonly arise during implementation. Numerical stability issues can emerge during the integration of highly sensitive systems, necessitating careful choice of tolerance parameters and integration schemes. Mesh quality must be maintained throughout the computation, often requiring smoothing algorithms that preserve geometric accuracy. The physical printing process presents its own challenges, particularly in optimizing model orientation to minimize support structures while ensuring structural integrity. These challenges can be addressed through careful parameter selection and monitoring of intermediate results at each stage of the pipeline.

Success in reproducing these results relies heavily on understanding both the mathematical foundations and practical considerations of 3D printing. Users should familiarize themselves with basic concepts of dynamical systems theory and computer-aided manufacturing. The provided code repository contains detailed documentation, example scripts, and printable files for test cases. These can serve as starting points for new implementations. We encourage users to experiment with different systems and parameters, contributing their findings and improvements to the broader scientific community through our open-source platform.

\section{Code repository reference.}\label{sec:code}
All code implementing the algorithms and visualization pipeline described in this article is available in our open-source repository:


\begin{verbatim}
https://github.com/esander1789/3DPrintingInvariantManifolds
\end{verbatim}

The repository includes:
\begin{itemize}
    \item Wolfram Mathematica\textsuperscript{\textregistered} implementations of Algorithms~\ref{alg:loc1} and~\ref{alg:p2}.
    \item Example scripts for the Lorenz, Arneodo-Coullet-Tresser, and Langford systems.
    \item Utility functions for mesh/surface generation and file export.
    \item Documentation and usage examples.
    \item Sample parameter files and printable (stl/obj) files test cases.
\end{itemize}

Users are encouraged to submit issues, feature requests, and contribute improvements through pull requests. For detailed usage instructions and documentation, please refer to the repository's README file and accompanying documentation. The code is released under the MIT License to promote open scientific collaboration and reproducibility.

For questions, bug reports, or collaboration inquiries, please use the GitHub issue tracker or contact the authors directly.


\section*{Acknowledgements} E.F. was supported by the National Science Foundation under Grant No. DMS-2137947. E.S. was supported in part by the Simons Foundation under Award No. 636383. We are grateful to Dan Anderson who helped us with and took several of the photographs.

\section*{Author Contributions}

\begin{itemize}
\item E.S. and E.F.: Conceptualization; 
Formal analysis; 
Methodology; 
Visualization; 
Writing – original draft; 
Writing – review and editing. 

\item A.G. Oguedo Oliva and J.A. Seay: Methodology and 3D printing.
\item P. Bishop and S. Chenoweth: 3D printing.
\end{itemize}

\bibliographystyle{vancouver}
\bibliography{references.bib}

\begin{thebibliography}{10}

\bibitem{lucas20}
Lucas SK, Sander E, Taalman L.
\newblock Modeling {Dynamical} {Systems} for {3D} {Printing}.
\newblock Notices of the American Mathematical Society. 2020 Dec;67(11).
\newblock Available from:
  \url{https://www.ams.org/notices/202011/rnoti-p1692.pdf}.

\bibitem{bertacchini23}
Bertacchini F, Pantano PS, Bilotta E.
\newblock Jewels from chaos: {A} fascinating journey from abstract forms to
  physical objects.
\newblock Chaos: An Interdisciplinary Journal of Nonlinear Science. 2023
  Jan;33(1):013132.
\newblock Available from:
  \url{https://pubs.aip.org/cha/article/33/1/013132/2877613/Jewels-from-chaos-A-fascinating-journey-from}.

\bibitem{gagliardo20183d}
Gagliardo M.
\newblock 3d printing chaos.
\newblock In: Proceedings of Bridges 2018: Mathematics, Art, Music,
  Architecture, Education, Culture; 2018. p. 491-4.

\bibitem{osinga2004crocheting}
Osinga HM, Krauskopf B.
\newblock Crocheting the Lorenz manifold.
\newblock Mathematical Intelligencer. 2004;26(4):25-37.

\bibitem{meiss07}
Meiss JD.
\newblock Differential {Dynamical} {Systems}.
\newblock Society for Industrial and Applied Mathematics; 2007.
\newblock Available from:
  \url{http://epubs.siam.org/doi/book/10.1137/1.9780898718232}.

\bibitem{Guck}
Guckenheimer J, Holmes P.
\newblock Nonlinear oscillations, dynamical systems, and bifurcations of vector
  fields. vol.~42.
\newblock Springer Science \& Business Media; 2013.

\bibitem{Cabre2003a}
Cabr\'{e} X, Fontich E, de~la Llave R.
\newblock The parameterization method for invariant manifolds I: manifolds
  associated to non-resonant subspaces.
\newblock Indiana University Mathematics Journal. 2003;52(2):283-328.

\bibitem{fleurantin2020resonant}
Fleurantin E, James JM.
\newblock Resonant tori, transport barriers, and chaos in a vector field with a
  Neimark--Sacker bifurcation.
\newblock Communications in Nonlinear Science and Numerical Simulation.
  2020;85:105226.

\bibitem{Krauskopf2007}
Krauskopf B, Osinga HM, Galán-Vioque J, editors.
\newblock Numerical continuation methods for dynamical systems.
\newblock Springer; 2007.

\bibitem{dellnitz2016computation}
Dellnitz M, Hessel-Von~Molo M, Ziessler A.
\newblock On the computation of attractors for delay differential equations.
\newblock Journal of Computational Dynamics. 2016;3(1):93-112.

\bibitem{Krauskopf05}
Krauskopf B, Osinga HM, Doedel EJ, Henderson ME, Guckenheimer J, Vladimirsky A,
  et~al.
\newblock A survey of methods for computing (un)stable manifolds of vector
  fields.
\newblock Internat J Bifur Chaos Appl Sci Engrg. 2005;15(3):763-91.
\newblock Available from: \url{https://doi.org/10.1142/S0218127405012533}.

\bibitem{henderson2005computing}
Henderson ME.
\newblock Computing invariant manifolds by integrating fat trajectories.
\newblock SIAM Journal on Applied Dynamical Systems. 2005;4(4):832-82.

\bibitem{haro16}
Haro {\`A}, Canadell M, Figueras JL, Luque A, Mondelo JM.
\newblock The {Parameterization} {Method} for {Invariant} {Manifolds}. vol. 195
  of Applied {Mathematical} {Sciences}.
\newblock Cham: Springer International Publishing; 2016.
\newblock Available from:
  \url{http://link.springer.com/10.1007/978-3-319-29662-3}.

\bibitem{shahrubudin2019overview}
Shahrubudin N, Lee TC, Ramlan R.
\newblock An overview on 3D printing technology: Technological, materials, and
  applications.
\newblock Procedia manufacturing. 2019;35:1286-96.

\bibitem{lee2017fundamentals}
Lee JY, An J, Chua CK.
\newblock Fundamentals and applications of 3D printing for novel materials.
\newblock Applied materials today. 2017;7:120-33.

\bibitem{segerman20123d}
Segerman H.
\newblock 3D printing for mathematical visualisation.
\newblock The Mathematical Intelligencer. 2012;34(4):56-62.

\bibitem{zastrow2020new}
Zastrow M.
\newblock The new 3D printing.
\newblock Nature. 2020;578(7793):20-3.

\bibitem{jandyal20223d}
Jandyal A, Chaturvedi I, Wazir I, Raina A, Haq MIU.
\newblock 3D printing--A review of processes, materials and applications in
  industry 4.0.
\newblock Sustainable Operations and Computers. 2022;3:33-42.

\end{thebibliography}

\vfill\eject

\end{document}